\documentclass[12pt,reqno]{amsart}
\usepackage{amsmath}
\usepackage{amssymb}
\usepackage{amscd}
\usepackage{mathrsfs}

\newtheorem{theorem}{Theorem}

\newtheorem{theoremb}{Theorem}
\newtheorem{theoremc}{Theorem}

\newtheorem{theoreme}{Theorem}

\newtheorem{rk}[theoremc]{Remark}
\newtheorem{cor}[theoremb]{Corollary}
\newtheorem{lem}[theoreme]{Lemma}

\newcommand\bib[1]{\bibitem[#1]{#1}}

\renewcommand\a{\alpha}
\renewcommand\b{\beta}

\newcommand\C{{\mathcal C}}
\newcommand\CC{{\mathbb C}}
\newcommand\Cc{{\let\mathcal\mathscr\mathcal C}}
\newcommand\com[1]{}

\newcommand\D{{\mathcal D}}
\newcommand\Dd{{\let\mathcal\mathscr\mathcal D}}

\newcommand\E{\mathcal{E}}

\newcommand\g{{\frak g}}

\newcommand\h{{\frak h}}

\newcommand\La{\Lambda}
\newcommand\Ll{{\let\mathcal\mathscr\mathcal L}}
\newcommand\m{\frak m}

\newcommand\oo{\omega}
\newcommand\op[1]{\mathop{\rm #1}\nolimits}
\newcommand\ot{\otimes}
\newcommand\p{\partial}

\newcommand\R{{\mathbb R}}
\renewcommand\t{\tau}
\newcommand\sym{\op{sym}}

\newcommand\z{\sigma}
\newcommand\Z{{\mathbb Z}}

\begin{document}

 \title{Finite-dimensionality in Tanaka theory}
 \author{Boris Kruglikov}
 \address{Department of Mathematics and Statistics, University of Troms\o, Troms\o, 9037, Norway.
 \qquad\hfill\it{boris.kruglikov@uit.no}}

 \begin{abstract}
In this paper we extend the Tanaka finiteness theorem and inequality for the
number of symmetries to arbitrary distributions
(differential systems) and provide several applications.
 \end{abstract}

 % MSC: 58A30, 34C14, 35A30, 58J72; 58J70, 34H05
 % thirdly: 58A17, 58A20, 34A05,

 \maketitle

%%%%%%%%%%%%%%%
\section*{Introduction and main results}

Let $\Delta$ be a regular completely non-holonomic distribution on a connected smooth manifold $M$,
i.e. a vector subbundle of $TM$. This paper concerns the Lie algebra $\op{sym}(\Delta)$ of its symmetries.

In \cite{T} N.\,Tanaka introduced a graded nilpotent Lie algebra (GNLA) $\m_x$ and its
algebraic prolongation $\g_x=\hat\m_x$ at every point $x\in M$ to majorize the symmetry
algebra of $\Delta$ (precise definitions follow in the next section).
One of the purposes of this paper is to prove

 \begin{theorem}\label{Thm1}
The Lie algebra $\op{sym}(\Delta)$ of the symmetries of $\Delta$ satisfies:
 $$
\dim\sym(\Delta)\le\op{sup}_M\dim\g_x.
 $$
 \end{theorem}

This statement was proven in \cite{T} (Corollary of Theorem 8.4) for strongly regular systems,
i.e. such $\Delta$ that the GNLA $\m_x$ does not depend on $x$.
Generic distributions on high-dimensional manifolds
(for instance, for rank 2 distributions starting from dimension 8) fail to satisfy this property,
but we relax the assumption to usual regularity.

 \begin{rk}\label{rk1}
In fact, the equality in the above theorem, provided that all $\g_x$ are finite-dimensional
Lie algebras, is attained only in one case: when the distribution $\Delta$ is strongly regular
and flat in the Tanaka sense (then it is locally isomorphic to the standard model of \cite{T}).
 \end{rk}

 \begin{rk}\label{rk2}
If we assume $\dim\g_x<\infty$ $\forall x\in M$, then the above inequality refines to
 \begin{equation}\label{ineq-rk2}
\dim\sym(\Delta)\le\op{inf}_M\dim\g_x.
 \end{equation}
 \end{rk}

Our approach here is to elaborate upon the Tanaka theory of symmetries for distributions
and relate it to the Spencer theory of formal integrability for PDEs.
Existence of this relation is natural since both theories are abstract versions and generalizations
of the Cartan equivalence method \cite{St,SS,Y}.

 \begin{cor}\label{cor}
Suppose for every $x\in M$ there exists no grading 0 derivation of the GNLA $\m_x^\CC=\oplus_{i<0}\g_i^\CC$,
which has rank 1 in $\op{gl}(\g_{-1}^\CC)$ and acts trivially on $\g_i^\CC$, $i<-1$.
Then the Lie algebra $\sym(\Delta)$ is finite-dimensional.
 \end{cor}

The meaning of the condition in the theorem is that the complex characteristic variety of
(the prolongation-projection of the Lie equation corresponding to) $\Delta$ is empty.
When the distribution is strongly regular this assertion follows from \cite{T} (Corollary 2 of Theorem 11.1).

This latter theorem and corollary concern finite-dimensionality of the Lie algebra
$\g_x$ and are purely algebraic. That is the reason it is "if and only if" statement.
For just regular distributions (even strongly regular but non-flat) there are other
reasons for decrease of the size of the symmetry algebra (for instance, abundance
of the independent components of the curvature).

Our theorems imply a-priory knowledge of finite-dimensionality and size estimate for the algebra
$\op{sym}(\Delta)$. Here is one output (the derived distribution is defined in Section \ref{S1}).

 \begin{theorem}\label{Thm2}
Consider a distribution $\Delta$ of rank $n$. Assume that rank of the derived distribution
$\Delta_2$ is greater than $\frac{n(n-1)}2+2$ in the case $n>2$. For $n=2$ we assume that rank of
the 2nd derived distribution $\Delta_3$ is $5$. 
Then the Lie algebra $\op{sym}(\Delta)$ is finite-dimensional.
 \end{theorem}

Another application of Corollary \ref{cor} is investigation of distributions
with infinite algebra of symmetries. In Section \ref{S8} we construct a model that,
together with the operation of prolongation, allows to classify
all rank 2 and 3 distributions with infinitely many symmetries.

\medskip

\textsc{Acknowledgment.}
I am grateful to Ian Anderson, Keizo Yamaguchi and Igor Zelenko for useful discussions
during the mini-workshop organized by Ian Anderson in Utah State University in November 2009.

%%%%%%%%%%%%%%%
\section{Review of Tanaka theory}\label{S1}

Given a distribution $\Delta\subset TM$, its \emph{weak derived flag}
$\{\Delta_i\}_{i>0}$ is given via the module of its sections by
$\Gamma(\Delta_{i+1})=[\Gamma(\Delta),\Gamma(\Delta_i)]$ with $\Delta_1=\Delta$.

We will assume throughout this paper that our distribution is \emph{completely non-holonomic},
i.e. $\Delta_i=TM$ for $i\ge\kappa$, and we also assume that the flag of $\Delta$ is \emph{regular},
so that the ranks of $\Delta_i$ are constant (whence $\kappa$ is constant as well).

The quotient sheaf $\m=\oplus_{i<0}\g_i$, $\g_i=\Delta_{-i}/\Delta_{-i-1}$ (we let $\Delta_0=0$),
has a natural structure of graded nilpotent Lie algebra at
any point $x\in M$. The bracket on $\m$ is induced by the commutator
of vector fields on $M$. $\Delta$ is called \emph{strongly regular} if the GNLA $\m=\m_x$ does not depend
on the point $x\in M$.

The \emph{growth vector} of $\Delta=\g_{-1}$ is the sequence of dimensions\footnote{In other sources, it is $(\dim\Delta_1,\dim\Delta_2,\dots)$.}\linebreak
 $(\dim \g_{-1},\dim \g_{-2}\dots)$ depending on $x\in M$.

The Tanaka prolongation $\g=\hat\m$ is the graded Lie algebra with negative graded part $\m$
and non-negative part defined successively by
 $$
\g_k=\{u\in\bigoplus\limits_{i<0}\g_{k+i}\ot\g_i^*:
u([X,Y])=[u(X),Y]+[X,u(Y)],\ X,Y\in\m\}.
 $$
Since $\Delta$ is bracket-generating, the algebra $\m$ is fundamental, i.e.
$\g_{-1}$ generates the whole GNLA $\m$, and therefore the grading $k$ homomorphism $u$
is uniquely determined by the restriction $u:\g_{-1}\to\g_{k-1}$.

The space $\g=\oplus\g_i$ is naturally a graded Lie algebra, called
the {\em Tanaka algebra\/} of $\Delta$, and to indicate dependence on the point $x\in M$, we
will write $\g=\g_x$ (also the value of a vector field $Y$ at the point $x$ will be denoted by $Y_x$).

Alternatively the above symbolic prolongation can be defined via Lie algebra cohomology with coefficients:
$\g_0=H^1_0(\m,\m)$, $\g_1=H^1_1(\m,\m\oplus\g_0)$ etc, where the subscript indicates
the grading \cite{AK}. The prolongation of $\m$ is $\g=\m\oplus\g_0\oplus\g_1\oplus\dots$.
In particular, if $\g$ is finite-dimensional, then $H^1(\m,\g)=0$ in accordance with \cite{Y}.

In addition to introducing the Lie algebra $\g$, which majorizes the symmetry algebra of $\Delta$,
the paper \cite{T} contains the construction of an important ingredient to the equivalence problem
-- an absolute parallelism on the prolongation manifold of the structure.

Distribution is locally flat if the structure functions of the absolute parallelism vanish.
Then the distribution $\Delta$ is locally diffeomorphic
to the standard model on the Lie group corresponding to $\m$, see \cite{T}.

%%%%%%%%%%%%%%%
\section{Lie algebra sheaf of a distribution}\label{S2}

Instead of considering the global object (vector fields) $\op{sym}(\Delta)\subset\Dd(M)$
let us study the more general structure of local symmetries. Namely we consider the
Lie algebras sheaf (LAS) $\Ll$ of germs of vector fields preserving the distribution $\Delta$.
Regularity of the latter implies that $\Ll$ is a sub-sheaf of $\Dd_\text{loc}(M)=\Gamma_\text{loc}(TM)$.

Let $\Ll(x)$ be its stalk at $x\in M$ and $\Ll(x)^0$ the isotropy subalgebra.
Note that in terms of the evaluation map $\op{ev}_x:\Ll(x)\to T_xM$ the latter is $\op{ev}_x^{-1}(0)$.
The LAS $\Ll$ is transitive if $\op{ev}_x$ is onto and in this case all stalks $\Ll(x)$ are isomorphic.
In the opposite case (provided $\Ll$ is integrable and regular) the Frobenius theorem implies that
$M$ is foliated by the leaves of $\op{ev}_x(\Ll(x))$; the stalks are constant along it.

We define the dimension of LAS to be
 $$
\dim\Ll=\op{sup}_M(\op{rank}[\op{ev}_x]+\dim\Ll(x)^0).
 $$
Since globalization can only decrease objects (see the Appendix for precise conditions), we will have:
 $$
\op{dim}\op{sym}(\Delta)\le\dim\Ll.
 $$

There are two decreasing filtrations of the stalk $\Ll(x)$ of the LAS $\Ll$ for every $x\in M$.

The first filtration is defined in the transitive case
by the rule \cite{SS}: $\Ll(x)_*^i=\Ll(x)$ for $i<0$, $\Ll(x)_*^0=\Ll(x)^0$
and
 $$
\Ll(x)_*^{i+1}=\{X\in\Ll(x)_*^i:[X,\Ll(x)]\subset\Ll(x)_*^i\}\text{ for }i\ge0.
 $$
This definition works in the abstract setting and in general there are existence and realization theorems
\cite{SS}. In our case this filtration is closely related to the jets. In fact,
$\Ll(x)_*^i=\Ll(x)\cap\mu_x^{i+1}\cdot\Dd_\text{loc}(M)$, where $\mu_x\subset C^\infty(M)$ is the maximal ideal.
This latter definition of $\Ll(x)_*^i$ ($i\ge0$) works nicely in the intransitive case as well.

The second filtration in the transitive case is defined as follows \cite{T}:
$\Ll(x)^i=\op{ev}_x^{-1}(\Delta_{-i})$ for $i<0$,
$\Ll(x)^0$ as above and
 $$
\Ll(x)^{i+1}=\{X\in\Ll(x)^i:[X,\Ll(x)^{-1}]\subset\Ll(x)^i\}\text{ for }i\ge0.
 $$
In particular, $\Ll(x)^{-i}=\Ll(x)$ for $i\ge\kappa$.

The two filtrations are related by the following diagram:
 $$
 \begin{array}{cccccccc}
\Ll(x) \!&= \!&\dots \!&= \Ll(x)^{-1}_* \!&\supset \Ll(x)^0_* \!&\supset \Ll(x)^1_* \!&\supset \Ll(x)^2_* \!&\supset \dots\\
|| & & & \ \cup \!\!& \ \ || \!\!& \ \ \cap \!\!& \ \ \cap \!\!& \\
\Ll(x)^{-\kappa} \!&\supset \!&\dots \!&\supset \Ll(x)^{-1} \!&\supset \Ll(x)^0 \!&\supset \Ll(x)^1 \!&\supset \Ll(x)^2 \!&\supset \dots
 \end{array}
 $$

In the intransitive case the situation with the second filtration is more complicated and can be resolved as follows.
The non-positive terms keep the same value.

For $X\in\Ll(x)^0$ and $Y\in\Dd(M)$ the value $\Psi_X^1(Y)=[X,Y]_x\in T_xM$ depends on $Y_x\in T_xM$ only.
Thus $\Psi_X^1:T_xM\to T_xM$ is a linear morphism and it maps $\Delta_i$ to $\Delta_i$, $i>0$.
We define
 $$
\Ll(x)^1=\{X\in\Ll(x)^0:\Psi_X^1|_{\Delta_x}=0\}.
 $$
Thus $X\in\Ll(x)^1$ are characterized by the property $\Psi_X^1:\Delta_i\to\Delta_{i-1}$.

For $X\in\Ll(x)^1$ and $Y,Z\in\Dd(M)$ the value $\Psi_X^2(Y,Z)=[[X,Y],Z]_x\in T_xM$ depends on
$Y_x,Z_x\in T_xM$ provided $Y_x,Z_x\in\Delta_x$.

Indeed, since $\R$-linearity of $\Psi_X^2$ is obvious, consider multiplication by a function $f\in C^\infty(M)$:
 \begin{multline*}
[[X,fY],Z]=f[[X,Y],Z]+X(f)[Y,Z]\\
-Z(f)[X,Y]-X(Z(f))Y+[X,Z](f)Y.
 \end{multline*}
At the point $x$ all the terms except the first in the r.h.s.\ vanish, implying the claim. The claim
for the second argument is even easier, as the last term below vanishes at $x$:
 $$
[[X,Y],fZ]=f[[X,Y],Z]+[X,Y](f)Z.
 $$
Thus we have a linear morphism $\Psi_X^2:\Delta_x\otimes\Delta_x\to\Delta_x$,
which we can also treat as a map $\Delta_x\to\Delta_x^*\otimes\Delta_x$,
$Z\mapsto\Psi_X^2(\cdot,Z)$.

We define
 $$
\Ll(x)^2=\{X\in\Ll(x)^1:\Psi_X^2|_{\Delta_x\ot\Delta_x}=0\}.
 $$
The Jacobi identity implies that $\Psi^1_X|_{\Delta_2}=0$ for $X\in\Ll(x)^2$
and more generally $\Psi_X^1:\Delta_i\to\Delta_{i-2}$ for such $X$ and $i\ge2$.

Continuing in the same way we define $\Ll(x)^i$ for $i>0$ and we inductively get
the multi-linear maps $\Psi_X^{i+1}:\ot^{i+1}\Delta_x\to T_xM$, $X\in\Ll(x)^i$ by the formula
 $$
\Psi_X^{i+1}(Y_1,\dots,Y_{i+1})=[[..[[X,Y_1],Y_2],..],Y_{i+1}]_x.
 $$
In the next section we extend $\Psi_X^{i+1}$ to various arguments and show
that $\Psi_X^{i+1}:\ot^{i+1}\Delta_x\to\Delta_x$.
Then the next filtration term is
 $$
\Ll(x)^{i+1}=\{X\in\Ll(x)^i:\Psi_X^{i+1}|_{\Delta_x\ot..\ot\Delta_x}=0\}.
 $$

These two filtrations define the same topology on $\Ll(x)$ because
 \begin{equation}\label{2filtr}
\Ll(x)^{i\cdot\kappa}\subset\Ll(x)_*^i\subset\Ll(x)^i,\ \ i\ge0.
 \end{equation}

%%%%%%%%%%%%%%%
\section{Investigation of the second filtration}\label{S3}

We need first to prove the following characterization of $X\in\Ll(x)^i$, $i\ge0$.
Let $Z_{jl}$ be a frame around $x$, compatible with the weak flag, i.e. $Z_{1l}(x)$ is a basis of
$\g_{-1}=\Delta_x$, $Z_{2l}(x)$ produces a basis of $\g_{-2}=\Delta_2/\Delta_1$ etc.
Decompose a local symmetry $X=\sum f_{jl}\,Z_{jl}$.

Let us say that function $f\in C^\infty(M)$ belongs to $\mu_{\Delta,x}^{k+1}$ if
$Y_1\cdots Y_t(f)$ vanishes at $x$ for all $Y_j\in\Gamma(\Delta)$, $t\le k$
(for $k=0$ we have: $\mu_{\Delta,x}=\mu_x$).

 \begin{lem}\label{L0}
If for $i\ge0$ $X\in\Ll(x)^i$, then $f_{jl}\in\mu_{\Delta,x}^{i+j}$, $1\le j\le\kappa$.
 \end{lem}

 \begin{proof}
$X$ is a local symmetry if for any $Y_{i_1}\in\Gamma(\Delta)$ the Lie bracket
 $$
[X,Y_{i_1}]=\sum(f_{jl}[Z_{jl},Y_{i_1}]-Y_{i_1}(f_{jl})Z_{jl})
 $$
is also a section of $\Delta$. This equality $\op{mod}\Delta_{j-1}$ implies the following decomposition
 \begin{equation}\label{fYi1}
Y_{i_1}(f_{jl})=\sum_{s\ge j-1}\a_{i_1jl}^{st}f_{st},\quad
j>1.
 \end{equation}
Differentiation of this implies
 \begin{equation}\label{fYi2}
Y_{i_1}Y_{i_2}(f_{jl})=\sum_{s\ge j-2}\a_{i_1i_2jl}^{st}f_{st},\quad j>2;
 \end{equation}
 \begin{equation}\label{fYi3}
Y_{i_1}Y_{i_2}(f_{2l})=\sum_{s,t}\a_{i_1i_22l}^{st}f_{st}+\b_{i_1i_2l}^tY_{i_1}(f_{1t}).
 \end{equation}
Proceeding we get generally for $j>1$
 \begin{equation}\label{fYi5}
Y_{i_1}\cdots Y_{i_r}(f_{jl})=\sum_{s\ge j-r}\a_{i_1\dots i_rjl}^{st}\,f_{st}+
\sum_{k=1}^{r-j+1}\b_{i_1\dots i_rjl}^{q_1\dots q_kt}\, Y_{q_1}\cdots Y_{q_k}(f_{1t}),
 \end{equation}
where we assume $\sum_A^B=0$ if $A>B$; we can also assume that $(q_1\dots q_k)$ is an ordered subset of
$(i_1\dots i_r)$ though it plays no role in what follows.

Consider at first the case $i=0$. Then obviously $f_{jl}(x)=0$, i.e. $f_{jl}\in\mu^1_x$.
Substituting this into (\ref{fYi1}) we get $f_{jl}\in\mu^2_{\Delta,x}$ for $j\ge2$.
Next substitution into (\ref{fYi2}) yields $f_{jl}\in\mu^3_{\Delta,x}$ for $j\ge3$,
and continuing with (\ref{fYi5}) we get $f_{jl}\in\mu^j_{\Delta,x}$ for all $j$.

Now let us look at $i=1$. By definition $f_{1l}\in\mu^2_{\Delta,x}$. Then
(\ref{fYi3}) implies $f_{2l}\in\mu^3_{\Delta,x}$. Continuing with (\ref{fYi5}) we get
$f_{jl}\in\mu^{j+1}_{\Delta,x}$.

Now the pattern is clear and the general claim
$f_{jl}\in\mu_{\Delta,x}^{i+j}$ is easily obtained by induction.
 \end{proof}

 \begin{lem}\label{L1/2}
For $X\in\Ll(x)^i$ and $Y_1,\dots,Y_{i+1}\in\Delta_x$ we have:
 $$\Psi^{i+1}_X(Y_1,\dots,Y_{i+1})\in\Delta_x.$$
 \end{lem}

 \begin{proof}
Indeed Lemma \ref{L0} implies that for $X\in\Ll(x)^i$ and $Y_j\in\Gamma(\Delta)$
 $$
[[..[X,Y_1],\dots Y_i],Y_{i+1}]=\sum\pm Y_{i+1}\cdots Y_1(f_{1l})Z_{1l}+\dots,
 $$
where the omitted terms that are linear combinations of $f_{jl}$, $Y_{q_1}(f_{jl})$, \dots,
$Y_{q_1}\cdots Y_{q_i}(f_{jl})$ and $Y_{q_1}\cdots Y_{q_{i+1}}(f_{tl})$ for $t>1$,
all of which vanish at $x$. Thus the result evaluated at $x$ is a vector from $\Delta_x$.
 \end{proof}

Recall that $\Delta_s=0$ for $s\le0$ and $\Delta_s=TM$ for $s\ge\kappa$.
The formulae of the previous section for any $X\in\Ll(x)$ yield the maps ($j,s_\nu>0$)
 $$
\Psi_X^j:\Gamma(\Delta_{s_1})\ot\dots\ot\Gamma(\Delta_{s_j})
\to\Gamma(\Delta_{s_1+...+s_j}).
 $$

 \begin{lem}\label{L1}
If $X\in\Ll(x)^i$ and $Y_t\in\Gamma(\Delta_{s_t})$ for $1\le t\le j\le i+1$, then
 $$
\Psi_X^j(Y_1,\dots,Y_j)_x\in\Delta_{s_1+\dots+s_j-i}.
 $$
 \end{lem}

 \begin{proof}
For $i=j$ and $s_1=\dots=s_j=1$ this is just the definition as $\Delta_0=0$.
Without loss of generality we can take decomposable
$Y_t=[[..[V_{t1},V_{t2}],..],V_{ts_t}]_x$, $V_{t\a}\in\Gamma(\Delta)$, $1\le t\le j$,
$1\le\a\le s_t$. Then
as a result of the Jacobi identity $\Psi^j_X(Y_1,\dots,Y_j)$ decomposes into a linear combination
of the terms
 $$
[[..[[X,V_{\a_1}],V_{\a_2}]..],V_{\a_r}],\quad r=s_1+\dots+s_j.
 $$
For $r\le i$ this vanishes at $x$, i.e. the result belongs to $\Delta_0$.

For $r>i$ we deduce from Lemma \ref{L0} for $X=\sum f_{jl}Z_{jl}$:
 $$
[[..[[X,V_{\a_1}],V_{\a_2}]..],V_{\a_r}]_x=\sum_{j\le r-i}\pm V_{\a_r}\cdots V_{\a_1}(f_{jl}Z_{jl})_x\in
\Delta_{r-i}.
 $$

Finally for $i=j-1$ the claim follows from Lemma \ref{L1/2} and the extension argument of this proof.
 \end{proof}

For $X\in\Ll(x)^i$ the multi-linear map $\Psi_X^j(Y_1,\dots,Y_j)$ of Lemma \ref{L1} is not $C^\infty(M)$-linear in $Y_t$,
% but when we quotient the values to $\g_{i-s_1-...-s_j}=\Delta_{s_1+\dots+s_j-i}/\Delta_{s_1+\dots+s_j-i-1}$ it is.
but considered with the values modulo $\Delta_{s_1+\dots+s_j-i-1}$ it is.
Thus it induces the map (all spaces evaluated at $x$)
 $$
\Psi_X^j:\Delta_{s_1}\ot\dots\ot\Delta_{s_j}\to\g_{i-s_1-...-s_j}
 $$
with the kernel $\sum_{t_1+...+t_j<s_1+...+s_j}\!
\Delta_{t_1}\ot\dots\ot\Delta_{t_j}$, whence the map
 \begin{equation}\label{Psi}
\Psi_X^j:\g_{-s_1}\ot\dots\ot\g_{-s_j}\to\g_{i-s_1-...-s_j}.
 \end{equation}

 \begin{rk}
For $X\in\Ll(x)^i$ and $i\ge(j-1)\kappa$ the map
 $$
\Psi_X^j:\Delta_{s_1}\ot\dots\ot\Delta_{s_j}\to\Delta_{s_1+...+s_j-i},
 $$
is tensorial, so in this case we do not need to quotient.
 \end{rk}

 \begin{lem}\label{L2}
The second filtration respects the Lie brackets:
 $$
[\Ll(x)^i,\Ll(x)^j]\subset\Ll(x)^{i+j}.
 $$
 \end{lem}

 \begin{proof}
Let $X'\in\Ll(x)^i,X''\in\Ll(x)^j$. If $i,j\le0$, then $X'_x\in\Delta_i$, $X''_x\in\Delta_j$ and
so $[X',X'']_x\in\Delta_{i+j}$, i.e. $[X',X'']\in\Ll(x)^{i+j}$.

Consider $i,j\ge0$. Then for $Y_1,\dots,Y_{i+j}\in\Gamma(\Delta)$ we have with some factors $\nu_t$:
 \begin{multline*}
[[..[[X',X''],Y_1]..,Y_{i+j}]=\\
\sum_{t;\ \z\in S_{i+j}}\nu_t\,[..[..[X',Y_{\z(1)}]..,Y_{\z(t)}],[..[X'',Y_{\z(t+1)}]..,Y_{\z(i+j)}]]
 \end{multline*}
If a term in the above summation has $t\ge i$, then $i+j-t\le j$ and for
$Z=[..[X',Y_{\z(1)}]..,Y_{\z(t)}]$ with $Z_x\in\Delta_{t-i}$ Lemma \ref{L1} implies:
$[[..[X'',Y_{\z(t+1)}]..,Y_{\z(i+j)}],Z]_x=0$.
If $t<i$, then a symmetric ($X'\leftrightarrows X''$) argument applies.

Let finally $i>0,j<0$, $i+j\ge0$. Then
 \begin{multline*}
[[..[[X',X''],Y_1]..,Y_{i+j}]=\\
\sum_{t,\z\in S_{i+j}}\nu_t\,[..[X',Y_{\z(1)}]..,Y_{\z(t)}],Z(X'',Y_{\z(t+1)},..,Y_{\z(i+j)})]]
 \end{multline*}
where $Z(..)=[..[X'',Y_{\z(t+1)}]..,Y_{\z(i+j)}]$ has $Z_x\in\Delta_{i-t}$, so that the result follows from
Lemma \ref{L1}.

If $i+j<0$, then the claim follows from Lemma \ref{L0}.
 \end{proof}

%%%%%%%%%%%%%%%
\section{Formal Lie algebra of symmetries}\label{S4}

In Section \ref{S2} we introduced two compatible filtrations (\ref{2filtr}).
The corresponding \textit{formal Lie algebra} is
 $$
L^\Delta_x=\lim_{i\to+\infty}\Ll(x)/\Ll(x)^i=\lim_{i\to+\infty}\Ll(x)/\Ll(x)_*^i.
 $$
This algebra has two gradations (jets and weighted jets\footnote{Weighted jets play a
crucial role in Morimoto's approach to the equivalence problem \cite{M}.}) corresponding to
the above two filtrations. In fact, the first grading is
 $$
\op{gr}_*(L^\Delta_x)=\oplus\,\bar g_i(x),\quad
\bar g_i(x)=\Ll(x)_*^{i-1}/\Ll(x)_*^i\ (i\ge0).
 $$
(the reason for the shift of indices will be clear in the next section).

The second grading is the following
 $$
\op{gr}(L^\Delta_x)=\oplus\,\mathbf{g}_i(x),\quad
\mathbf{g}_i(x)=\Ll(x)^i/\Ll(x)^{i+1}\ (i\ge-\kappa).
 $$
Notice the difference in the range of indices.
Both gradings have the induced Lie bracket, so that we have two graded Lie structures (which
might be different as Lie algebras from $L^\Delta_x$).

We clearly have the inclusion of GNLAs $\oplus_{i<0}\mathbf{g}_i\subset \m_x$, and now would like to
elaborate upon (\ref{Psi}) to understand the symbols $\mathbf{g}_i$, $i\ge0$.

The Tanaka symbol space $\g_0\subset\sum_{i<0}\g_i^*\ot\g_i$ is uniquely determined by
its restriction to $\g_{-1}$ and thus can be identified with its image
$\g_0\hookrightarrow\g_{-1}^*\ot\g_{-1}$.

Similarly, $\g_1\subset(\sum_{j<0}\g_{-1}^*\ot\g_j^*\ot\g_j)\oplus(\sum_{i<-1}\g_{i-1}^*\ot\g_i)$
can be identified with its image
$\g_1\hookrightarrow\g_{-1}^*\ot\g_0\hookrightarrow\g_{-1}^*\ot\g_{-1}^*\ot\g_{-1}$.

In the general case, $\g_i\hookrightarrow\ot^{i+1}\g_{-1}^*\ot\g_{-1}$ is a
monomorphism and we identify the symbol $\g_i$ ($i\ge0$) with its image.

Now notice that by (\ref{Psi}) every
$X\in\Ll(x)^i$ induces the linear map $\Psi_X^j:\ot^j\g_{-1}\to\g_{i-j}$.

 \begin{lem}\label{L3}
For $X\in\Ll(x)^i$ the element $\Psi_X^{i+1}\in\ot^{i+1}\g_{-1}^*\ot\g_{-1}$
belongs to (the image of) $\g_i$ ($i\ge0$).
 \end{lem}

 \begin{proof}
For $i=0$ the claim is obvious since $\Psi^1_X=\op{ad}_X\in\op{Der}_0(\m_x)$.

For $i=1$ the Jacobi identity implies the following symmetry of $\Psi_X^2$, with $Y,Z\in\Gamma(\Delta)$ so that
$[Y,Z]\in\Gamma(\Delta_2)$:
 $$
\Psi_X^2(Y,Z)-\Psi_X^2(Z,Y)=\Psi_X^1([Y,Z]).
 $$
Moreover this formula holds true if we understand $Y\in\g_j$, $Z\in\g_l$, $[Y,Z]\in\g_{j+l}$,
which provides us the extension $\Psi^2_X\in\sum_{i<0}\g_{i}^*\ot\g_{i+1}$
we seek. By the construction this element satisfies the Leibniz rule and so belongs to $\g_1$.
In other words, $\Psi^2_X|_{\Delta_x\ot\Delta_x}$
belongs to the image of $\g_1$ in $\g_{-1}^*\ot\g_{-1}^*\ot\g_{-1}$.

For the general $i>0$ the arguments are the same (induction), and the reason behind the
claim of Lemma \ref{L3} is that both the element $\Psi^{i+1}_X$ and the elements of $\g_i$
are constructed on the same principle using the same formula (cf. introduction of
the spaces $\g_i$ for $i\ge0$ in \cite{T,Y}).
 \end{proof}

Now we are ready for our main technical result.

 \begin{theorem}\label{hm3}
For all $i\in\Z$: $\mathbf{g}_i\subset\g_i$. Thus we get the monomorphism of the graded Lie
algebras:
 $$
\op{gr}(L^\Delta_x)=\oplus\mathbf{g}_i\hookrightarrow\oplus\g_i=\g_x
 $$
 \end{theorem}

 \begin{proof}
Consider the map
 $$
\Ll(x)^i\ni X\mapsto \Psi_X^{i+1}\in\g_i,\quad i\ge0.
 $$
Its kernel is $\Ll(x)^{i+1}$ and therefore we get the induced monomorphism $\mathbf{g}_i\to\g_i$.
Existence of this arrow for negative $i$ is obvious.

Since the Lie bracket in both cases is induced by the commutator of vector fields
and by Lemma \ref{L2} it respects the filtration, so it respects the gradation and the map is a homomorphism of Lie algebras.
 \end{proof}

Notice though that the formal Lie algebra $L^\Delta_x$ is not a Lie subalgebra of $\g_x$,
as can be seen by studying the sub-maximal examples in \cite{C$_1$}.

 \begin{cor}
At any point $x\in M$ we have $\dim L^\Delta_x\le\dim\g_x$.
 \end{cor}

Thus for the LAS $\Ll$ we achieved our claim. To prove it for
global/lo\-cal\footnote{This means over a small neighborhood, fixed for all vector fields.}
Lie algebra $\op{sym}(\Delta)$ we shall study the first grading too.

%%%%%%%%%%%%%%%
\section{Lie equation associated to a distribution}\label{S4}

We will use here jet-spaces and the geometric theory of PDE, for which we refer the reader to
\cite{Sp,KLV,KS,KL}. Consider the Lie equation $\E=\mathfrak{Lie}(K)$ of a geometric structure $K$.
For instance, if $K$ is a tensorial field, then this is the equation for its symmetries, i.e.
vector fields $X$ satisfying $L_X(K)=0$.

In the case of our current interest $K=\Delta$ and the Lie equation, considered as the
submanifold in jets, is
 $$
\E=\{j^1_xX:L_X(\Delta)_x\subset\Delta_x\,|\, x\in M\}\subset J^1(TM).
 $$

As is customary in the theory of overdetermined systems, finding compatibility conditions on the solutions
of $\E$ binds to applying the prolongation-projection method. Namely the prolongation
$\E^{(k-1)}\subset J^k$ is determined by the original equations and their derivatives up to order $k$.
It can happen that some projections $\pi_{k,l}:\E^{(k-1)}\to\E^{(l-1)}$ are not surjective
(i.e. there is a differential corollary of lower order),
then we take the image as the new system of equations and apply the prolongations and projections again.
The procedure is finite due to Cartan-Kuranishi theorem \cite{Ku} and the output is formally integrable,
meaning that it possesses a formal series solution through every regular point
(or local analytic solution for analytic $\E$).

Calculation of prolongation-projection in general is a difficult task, and $\mathfrak{Lie}(\Delta)$
is not an exception. Hopefully, we can guess an equation squeezed in between $\E$ and its $k$-th
derived $\pi_{k,1}(\E^{(k-1)})\subset J^1$. This is so because of the obvious fact that if $k$-jet of $X$
preserves $\Delta$ at $x$, then $j^1_xX$ preserves its weak derived flag up to order $k$.

Thus instead of studying symmetries of $\Delta$ we can equally well study symmetries of the derived
flag $\mathcal{W}_\Delta=\{\Delta_i\}_{i>0}$. Denote by $\E_\Delta$ the corresponding Lie equation
 $$
\mathfrak{Lie}(\mathcal{W}_\Delta)=
\{j^1_xX:L_X(\Delta_i)_x\subset(\Delta_i)_x\,\forall i>0\,|\,x\in M\}\subset J^1(TM).
 $$
Denote by $\bar\E$ the result of the prolongation-projection. This consists of the subsets
(submanifolds with singularities)
 $$
\bar\E_i=\lim_{j\to\infty}\pi_{j+1,i}(\E_\Delta^{(j)})\subset J^i(TM).
 $$

Symbols of this equation are the vector spaces (we use linearity of $\bar\E$
which simplifies the general formulae)
 $$
\bar g_i=\op{Ker}(\pi_{i,i-1}:\bar\E_i\to\bar\E_{i-1}).
 $$
Clearly these are subspaces of the symbols $g_i$ of the original Lie equation $\E_\Delta$:
 \begin{equation}\label{II-3}
\bar g_i\subset g_i\subset S^iT^*_xM\ot T_xM.
 \end{equation}
In particular, $\bar g_0$ is the tangent to the orbit of the symmetry group action
(the whole $T_xM$ in the transitive case).

Ultimately the infinite jets $\bar\E_\infty$ correspond to $L^\Delta_x$
(another descriptions is this: $\bar\E_\infty$ consists of formal vector fields 
preserving all differential invariants of $\Delta$).
We shall relate the Spencer symbols $\bar g_{i+1}$ to the Tanaka symbols $\g_i$, $i\ge0$.

Consider the decreasing filtration $\Ll(x)^i\cap\Ll(x)^j_*$ of the space $\Ll(x)^i$.
It produces the grading (the isomorphisms below are not natural and respect only the
linear structure)
 $$
\Ll(x)^i\simeq\bigoplus_{j\ge0}\frac{\Ll(x)^i\cap\Ll(x)^j_*}{\Ll(x)^i\cap\Ll(x)^{j+1}_*}=
\bigoplus_{j=0}^{i-1}\frac{\Ll(x)^i\cap\Ll(x)^j_*}{\Ll(x)^i\cap\Ll(x)^{j+1}_*}\oplus\bigoplus_{j\ge i}\bar g_{j+1}.
 $$
This implies for $i\ge0$

 \begin{equation}\label{BigO1}
\mathbf{g}_i=\Ll(x)^i/\Ll(x)^{i+1}\simeq\bigoplus_{j\le i}h_{ij},
 \end{equation}
where $h_{ij}=(\Ll(x)^i\cap\Ll(x)^j_*)/(\Ll(x)^i\cap\Ll(x)^{j+1}_*+\Ll(x)^{i+1}\cap\Ll(x)^j_*)$.
Note that $h_{ij}=0$ if either $i<j$ or $i\ge(j+1)\kappa$ ($j\ge0$).

Similarly the decreasing filtration $\Ll(x)^i\cap\Ll(x)^j_*$ of the space $\Ll(x)^j_*$ yields
 \begin{equation}\label{BigO2}
\bar g_{j+1}=\Ll(x)^j_*/\Ll(x)^{j+1}_*\simeq\bigoplus_{i\ge j}h_{ij}.
 \end{equation}

As a by-product of calculations in the previous sections,
we can interpret the vector space $h_{ij}$ as a subspace in
 $$
\sum_{s_1+{}\dots{}+s_{j+1}=i+t}\g_{-s_1}^*\ot\dots\ot\g_{-s_{j+1}}^*\ot\g_{-t}\quad(t,s_\nu\ge1).
 $$

 \begin{theorem}\label{hm4}
The Spencer symbols $\bar g_i$ of the equation $\bar\E$ are related to the
Tanaka symbols $\mathbf{g}_i$ via (\ref{BigO1})-(\ref{BigO2}).
This yields a (noncanonical) monomorphism of $\oplus_{j>0}\bar g_j$ into $\oplus_{i\ge0}\g_i$.
 \end{theorem}
The claim of the theorem follows from the inclusions $\mathbf{g}_i(x)\subset\g_i(x)$ of Theorem \ref{hm3}.
These are strict and the only case, when we have equalities for all $x$ is the Tanaka flat
distribution (for finite-dimensional $\g$)\footnote{In this case Theorem \ref{hm4} gives two gradings
on the space of symmetries of the standard model $\Delta=\Delta_{\m}$, but only $\op{gr}(L^\Delta_x)$
yields the  Lie algebra structure of $\op{sym}(\Delta)$.}.

Indeed, equalities everywhere mean that the Lie equation of the GNLA $\m_x$
is formally integrable, and the prolongation-projection does not decrease its symbols.
In particular, $\bar g_0(x)=T_xM$, the pseudogroup of symmetries is transitive, so that the distribution
is strongly regular and the result follows. This justifies Remark \ref{rk1}.

 \begin{cor}
For a regular distribution $\Delta$ (not necessarily strongly regular)
 $$
\sum_{i>0}\dim \bar g_i\le \sum_{j\ge0}\dim\g_j.
 $$
 \end{cor}

This implies that if the map $j_x^\infty:\Ll(x)^0\to J^\infty_x(TM)$ is injective
(we justify the assumption in the next section), then
 $$
\dim\Ll(x)\le\sum_{j=-\kappa}^\infty\dim\g_j(x),
 $$
and so $\dim\Ll\le\sup_M\sum_j\dim\g_j$.

%%%%%%%%%%%%%%%
\section{Proof of Theorem \ref{Thm1} and beyond}\label{S4.5}

We can suppose that $\dim\g_x$ is finite at every point $x\in M$, because else the inequality
in Theorem \ref{Thm1} is trivial.
 \begin{lem}
The function $x\mapsto\dim\g_x$ is upper semi-continuous.
 \end{lem}

 \begin{proof}
Indeed, $\g_x$ is obtained from $\m_x$ by certain linear algebra rules 
(in \cite{T,Y} the positive grades $\g_i$ are defined successively, but this can be easily
modified to obtain $\g_+$ via $\m$ at once).
Since ranks of matrices can only drop in the limit process, the result follows.
 \end{proof}
Thus $\dim\g_x$ attains a maximum in any compact domain $\bar U$.

Let $N$ be the number such that $\g_i(x)=0$ for all $i\ge N$ and $x\in\bar U$.
Then (\ref{BigO1})+(\ref{BigO2}) imply that $\bar g_i(x)=0$ for
all $i>N$ and $x\in\bar U$. In fact, the whole prolongation-projection process
is not required, but as the condition that $X$ preserves the GNLA $\m_x$ structure is obtained
from the original $\mathfrak{Lie}(\Delta)$ in a finite number of steps and
$\g_x$ is obtained via $\m_x$ by algebraic prolongation we conclude:
In finite number of steps of prolongation-projection the Lie equation becomes
of finite type at all points $x\in\bar U$.

Consequently by the results of Theorem \ref{T-app} from the Appendix,
there are no symmetries in $\bar U$ flat at some point.
Since $\bar U\subset M$ is arbitrary, we conclude the result for the whole $M$
and hence the map
 $$
j_x^\infty:\op{sym}(\Delta)\to\bar\E_x^\infty
 $$
is injective for every point $x\in M$.

Thus results of Section \ref{S4} imply the inequality
 $$
\dim\op{sym}(\Delta)\le\dim\Ll(x)\quad \forall x\in M
 $$
and consequently we prove inequality (\ref{ineq-rk2}) of Remark \ref{rk2}, which implies in
turn Theorem \ref{Thm1}.
\qed

 \begin{rk}
If the distribution is not strongly regular, then the leaves on $M$ through a typical point $x$
(obtained by fixing the invariants) have codimension $r=n-\op{rank}[\op{ev}_x]>0$ and we
refine inequality (\ref{ineq-rk2}) to
 \begin{equation*}
\dim\op{sym}(\Delta)\le\op{inf}_M\dim\g_x-r.
 \end{equation*}
 \end{rk}

\noindent{\bf Proof of Corollary \ref{cor}.}
We use Corollary 2 of Theorem 11.1 from \cite{T}, which states that if the subalgebra
 $$
\h_0=\{v\in g_0:[v,g_r]=0\ \forall r<-1\}
 $$
is of finite type as the subalgebra of $\op{gl}(\g_{-1})$, then $\g$ is finite-dimensional.

The subalgebra $\h_0\subset\op{gl}(\g_{-1})$ is of finite type iff its
complex characteristic variety is empty \cite{GQS}.

Thus all Tanaka algebras $\g_x$, $x\in M$, are finite-dimensional and the result follows.
\qed

\medskip

We can introduce the characteristic variety (the set of covectors satisfying the
defining relation is homogeneous and we projectivize it)
 $$
\op{Char}(\Delta)=\mathbb{P}\{p\in g_{-1}^*\setminus\{0\}:\exists q\in g_{-1}\setminus\{0\},\ p\ot q\in \mathfrak{h}_0\}
\subset \mathbb{P}\Delta^*.
 $$
Working over $\CC$ (this part of the theory is algebraic, and
complexification makes no problem) we obtain $\op{Char}^\CC(\Delta)\subset\mathbb{P}^\CC\Delta^*$.
The criterion of the theorem reformulates now as follows:
 $$
\op{Char}^\CC(\Delta)=\emptyset\ \Longrightarrow\ \dim\op{sym}(\Delta)<\infty.
 $$

If the set $\op{Char}^\CC(\Delta)$ is empty, then the complex characteristic variety of the differential
closure $\bar\E$ of the Lie equation $\E_\Delta$ is empty as well, but the reverse implication
is not generally true though it holds for a Tanaka flat distribution $\Delta$.

 \begin{rk}
One should be cautious as $\op{Char}^\CC(\bar\E)$ is usually smaller than
the initial complex characteristic variety $\op{Char}^\CC(\E_\Delta)$.
 \end{rk}

{\bf Example.} Consider the Tanaka flat rank 2 distribution in $\R^6$
of growth $(2,1,2,1)$ and maximal symmetry algebra $\mathfrak{p}_6$ of dimension 11.
In \cite{AK} (see also \cite{DZ}) it was shown that it corresponds to the Monge equation
 $$
\E_{1,3}:\ y'=(z''')^2,
 $$
namely $\Delta=\langle\D_x=\p_x+z_3^2\p_y+z_1\p_z+z_2\p_{z_1}+z_3\p_{z_2},\p_{z_3}\rangle$
in the standard jet-coordinates in mixed jets $J^{1,3}(\R,\R\times\R)\supset\E_{1,3}\simeq\R^6(x,y,z,z_1,z_2,z_3)$.

Moreover its algebra of symmetries $\op{sym}(\Delta)$ is isomorphic to the
Tanaka algebra $\g=\mathfrak{p}_6$ with $\g_0\simeq\R^3$, $\g_1\simeq\R^2$, $\g_2=0$,
and it has the following basis corresponding to elements of pure grade in $\g$:
 \begin{align*}
\g_{-4}:\quad & Z_0=\p_z,\\
\g_{-3}:\quad & Z_1=x\p_z+\p_{z_1},\ \ Y_0=\p_y,\\
\g_{-2}:\quad & Z_2=\tfrac{x^2}2\p_z+x\p_{z_1}+\p_{z_2},\\
\g_{-1}:\quad & Z_3=\tfrac{x^3}{3!}\p_z+\tfrac{x^2}2\p_{z_1}+x\p_{z_2}+\p_{z_3}+2z_2\p_y,\ \ S_0=\p_x\\
\g_0:\quad & Z_4=\tfrac{x^4}{4!}\p_z+\tfrac{x^3}{3!}\p_{z_1}+\tfrac{x^2}2\p_{z_2}+x\p_{z_3}+2(xz_2-z_1)\p_y,\\
           & S_1=x\p_x+\tfrac52z\p_z+\tfrac32z_1\p_{z_1}+\tfrac12z_2\p_{z_2}-\tfrac12z_3\p_{z_3},\\
           & R=y\p_y+\tfrac12z\p_z+\tfrac12z_1\p_{z_1}+\tfrac12z_2\p_{z_2}+\tfrac12z_3\p_{z_3},\\
\g_1:\quad & Z_5=\tfrac{x^5}{5!}\p_z+\tfrac{x^4}{4!}\p_{z_1}+\tfrac{x^3}{3!}\p_{z_2}+\tfrac{x^2}2\p_{z_3}+
2(\tfrac{x^2}2z_2-xz_1+z)\p_y,\\
           & S_2=x^2\p_x+9z_2^2\p_y+5xz\p_z+(5z+3xz_1)\p_{z_1}+\\
           & \qquad\qquad\qquad\qquad\qquad +(8z_1+xz_2)\p_{z_2}+(9z_2-xz_3)\p_{z_3}.
 \end{align*}

If we take the classes of these fields in $\Ll(x)^i_*/\Ll(x)^{i+1}_*$ we get:
 $$
\bar g_0=\langle[S_0],[Y_0],[Z_0],[Z_1],[Z_2],[Z_3]\rangle,\ \
\bar g_1=\langle[S_1],[S_2],[R],[Z_4],[Z_5]\rangle.
 $$
However this corresponds to the differential closure (via prolongation-projection) $\bar\E_\Delta$,
while the original Lie equation $\E_\Delta=\mathfrak{Lie}(\Delta)$ is bigger. In particular,
the element $[Z_5]=dz\ot\p_y\in\bar g_1\subset g_1$ is of rank 1. Further prolongation-projections
of $\E_\Delta$ yield $\bar g_2=0$ and kill this characteristic.

%%%%%%%%%%%%%%%
\section{Finite-dimensionality of the symmetry algebra}\label{S5}

Theorem \ref{Thm2} for $n=2$ was essentially established in \cite{AK}
(but in that paper we restricted to the strongly regular case)
by showing that $\h_0=0$ provided
the growth vector is $(2,1,2,\dots)$, and even in a more general case $(2,1,\dots,1,2,\dots)$.

Another proof is as follows. We can work over $\CC$ as this does not change the dimensions
of graded components. Suppose there is a nonzero element of rank 1
$\oo=p\ot\zeta\in\h_0\subset\g_{-1}^*\ot\g_{-1}$. Let $\xi$ be a complement to $\zeta$ in $\g_{-1}$.
Then $\oo$ acts trivially on $\g_{-2}=\R\cdot[\zeta,\xi]$ iff $p(\zeta)=0$. Furthermore
$\oo$ acts trivially on $\g_i$, $i<-2$, iff the operator
$\op{ad}_\zeta:\g_{i+1}\to\g_i$ is zero. Consequently $\dim\g_i=1$ for $i<-1$.

\smallskip

\noindent{\bf Proof of Theorem \ref{Thm2} for $n>2$.}
Again, working over $\CC$ and taking an element
of rank 1 $\oo=p\ot\zeta\in\h_0\subset\g_{-1}^*\ot\g_{-1}$ we observe that $\oo$ acts trivially
on $\g_{-2}$ iff $\op{ad}_\zeta|_{\op{Ann}(p)}=0$. As $\zeta$ can belong to $\op{Ann}(p)$,
this imposes $(n-2)$ restrictions.
Thus in the case the characteristic variety is non-empty, dimension of $\g_{-2}=\op{ad}(\La^2\g_{-1})$
does not exceed $\frac{n(n-1)}2-(n-2)=\frac{(n-1)(n-2)}2+1$.

Thus distributions with growth starting $(3,3,..)$, $(4,5,..)$, $(4,6,..)$ etc
have finite-dimensional symmetry algebras. \qed

\medskip

The a-priory knowledge of finite-dimensionality of the symmetry algebra can be enhanced
by estimation of its maximal size in many cases. For distributions of rank 2 it is done in \cite{DZ,AK}.

It can be also done for distributions $\Delta$ with free truncated GNLAs. Assume that
$\m_x$ is the free Lie algebra of step $k$, i.e. the only constraints are the Jacobi identity
and the requirement $\g_{-k-1}=0$. Then obviously $\g_0=\g_{-1}^*\ot\g_{-1}$. Moreover
by \cite{W} $\g_1=0$ provided $k>2$, $n>2$ or $k>3$, $n=2$. This implies
(\cite{W} concerned only homogeneous distributions, but the result holds for any $\Delta$
with truncated free $\m$ due to \cite{T} or Theorem \ref{Thm1}):
 $$
\dim\op{sym}(\Delta)\le \ell_n(k)+n^2
 $$
Here $\ell_n(k)=\dim\m_x$ can be calculated recursively or via the M\"obius function by \cite{Se}
 $$
k\,\ell_n(k) \,=\, n^k-\!\!\!\sum_{m|k,m<k}\!\! m\,\ell_n(m)
\,=\,\sum_{m|k}\mu(m)\,n^{k/m}.
 $$

For the exceptional cases we have:
$\dim\op{sym}(\Delta)\le 2n^2+n$, provided $k=2$, $n>2$
and the maximal symmetric case is given by $B_n=\mathfrak{o}(2n+1)$ \cite{Y,CN}.

For $k=3$, $n=2$ $\dim\op{sym}(\Delta)\le14$, the maximal dimension being realized only
in the case of $\op{Lie}(G_2)$ \cite{C$_1$}.
 %, while for $k=2$, $n=2$ $\op{sym}(\Delta)=\mathfrak{cont}(\R^3)$.

\smallskip

The estimate of Theorem \ref{Thm2} is sharp as there exist infinite-dimensional Lie algebras acting
as symmetries on distributions with growth vector beginning $(2,1,1,..),$ $(3,1,..)$, $(3,2,..)$, $(4,1,..)$,
$(4,2,..)$, $(4,3,..)$, $(4,4,..)$ etc, as discussed in the next Section.

%%%%%%%%%%%%%%%
\section{Infinite algebras of symmetries}\label{S8}

As in the rest of the paper we assume $\Delta$ totally non-holonomic. Notice however that
if the bracket-closure $\Delta_\infty\ne TM$, and $\mathcal{F}$ are the leaves of $\Delta_\infty$,
then existence of one symmetry transversal to $\mathcal{F}$ implies existence of
an infinite-dimensional space of such symmetries.

\subsection{Distributions of rank $n=2$}\label{SS81}

If a totally non-holonomic distribution with $\op{rank}(\Delta)=2$
has infinite-dimensional symmetry algebra then its growth vector starts $(2,1,1,\dots)$.

Recall that a Cauchy characteristic of a distribution is a symmetry tangent to this distribution.
The following statement is due to  E.\,Cartan \cite{C$_2$} (see \cite{AK} for another proof;
alternatively the vector field $\zeta$ is the one constructed in
the beginning of Section \ref{S5}).

 \begin{theorem}\label{ThDeProl}
The growth vector of a rank 2 distribution $\Delta$ is $(2,1,1,\dots)$
if and only if there exists a vector field $\zeta\in\Gamma(\Delta)$ which is a
Cauchy characteristic for the distribution $\Delta_2$.
In this, and only in this case, $\Delta$ is locally the prolongation of another
rank 2 distribution $\bar\Delta$.
 \end{theorem}

Here the (geometric) prolongation of a rank 2 distribution $\bar\Delta$ on $\bar M$ is
the projectivization $M=\mathbb{P}\bar\Delta=\{l_x\subset\bar\Delta_x:x\in\bar M\}$
with natural projection $\pi:M\to\bar M$ and the distribution $\Delta=\pi^{-1}_*(l_{\pi(x)})$.
Locally in $\bar M$ for $\bar\Delta=\langle U,V\rangle$ and $t\in S^1=\R^1/2\pi\Z$ we have:
$M=\bar M\times S^1$ and $\Delta=\langle \cos t\cdot U+\sin t\cdot V,\p_t\rangle$.

If the distribution $\bar\Delta$ has a preferred (vertical) section
$V\in\Gamma(\bar\Delta)$ [this means that the symmetries preserve it],
then the (affine) prolongation $\Delta=\bar\Delta^{(1)}$
is simply $\Delta=\langle U+t\,V,\p_t\rangle$ on $M=\bar M\times\R^1$, where $t$ is the coordinate on $\R^1$.
This coincides with (Spencer) prolongation in the geometric theory of PDEs (\cite{KLV,MZ}).

When $\Delta=\bar\Delta^{(1)}$, the operation $\Delta\mapsto\bar\Delta$ is called de-prolongation.
In the case of Theorem \ref{ThDeProl},
de-prolongation is the quotient of the first derived distribution by the Cauchy characteristic
 $$
(\bar M,\bar\Delta)=(M,\Delta_2)/\zeta.
 $$

The algebras $\sym(\bar\Delta)$ and $\sym(\bar\Delta^{(1)})$ are isomorphic
[indeed any symmetry on $\bar M$ induces the action on $t$ and thus lifts to $M$]. Therefore
Theorem \ref{Thm2} and a sequence of de-prolongations yield the following important statement
(it gives a-posteriori transitivity of the symmetry group action; with transitivity imposed
a-priori -- as an additional assumption -- the claim follows from Theorem 7.1 of \cite{MT}).

 \begin{theorem}
A regular germ of a rank 2 non-holonomic distribution $\Delta$ has infinite-dimensional symmetry
algebra if and only if $\Delta$ is equivalent to the Cartan distribution $\mathcal{C}_k$
on jet-space $J^k(\R,\R)$, $k=\dim M-2$.
 \end{theorem}

Recall that distributions $\Delta$ with strong growth vector $(2,1,1,\dots,1)$ are called
Goursat distributions. Their regular points can be characterized by the condition that
the growth vector is $(2,1,1,\dots,1)$ in both the weak and the strong sense.
Near such points the normal form
 $$
\mathcal{C}_k=\op{Ann}\{dy_i-y_{i+1}\,dx|0\le i<k\}\subset TJ^k(\R,\R)
 $$
is provided by the von Weber - Cartan theorem \cite{W,C$_2$}, see also \cite{GKR}.

\subsection{Distributions of higher rank}\label{SS82}

We indicate a local construction to produce a distribution with infinite-dimensional
(intransitive) symmetry algebra from any distribution.

For simplicity let's start with the case, when we extend a rank 2 distribution.
Locally any such distribution can be represented as the Cartan
distribution $\C_k$ for the Monge system, i.e. an underdetermined ODE $\E\subset J^k(\R,\R^m)$
given by $(m-1)$ equations.

Then we define $M=\E\times_\R J^l(\R,\R)\subset J^{k,l}(\R,\R^m\times\R)$ and
the distribution is $\Delta=(\C_k|_\E)\times_\R\C_l$.
In local coordinates if the equations in $\E$ are
$u^i_k=\psi^i(x,u^j,\dots,u^j_{k-1},v,\dots,v_k)$, $i,j=1,\dots,m-1$
($u^i=u^i(x)$, $v=u^m(x)$ are unknowns, $u^j_1=u^j_x$, $v_1=v_x$, $v_2=v_{xx}$ etc),
then $\C_k|_\E=\langle\D_x'=\p_x+u^i_1\p_{u^i}+\dots+\psi^i\p_{u^i_{k-1}}
+v_1\p_v+\dots+v_k\p_{v_{k-1}},\p_{v_k}\rangle$ and
$\C_l=\langle\D_x''=\p_x+w_1\p_w+w_2\p_{w_1}+\dots,\p_{w_l}\rangle$, so that
 $$
\Delta=\langle\D_x=\p_x+u^i_1\p_{u^i}+v_1\p_v+w_1\p_w+\dots,\p_{v_k},\p_{w_l}\rangle.
 $$
It is obvious that prolongation $\hat X$ of the vector field $X=f(w)\p_w$ is a symmetry
of $\Delta$ for any function $f(w)$.

Similarly a general regular rank $n$ distribution is realized locally as the Cartan distribution
$\C_k|_\E$ of an overdetermined PDE system $\E\subset J^k(\R^s,\R^m)$.
Then the distribution $\Delta=(\C_k|_\E)\times_{\R^s}\C_l$ on
$M=\E\times_{\R^s}J^l(\R^s,\R)\subset J^{k,l}(\R^s,\R^m\times\R)$ has an infinite algebra of
symmetries.

We can shrink the 2nd factor to a symmetric equation $\mathcal{R}\subset J^l(\R^s,\R)$,
for instance taking $\mathcal{R}=J^l(\R,\R)$ given by the equations $w_{x_2}=\dots=w_{x_s}=0$,
and still have infinitely many symmetries $\hat X$, $X=f(w)\p_w$ for the distribution
$\Delta=(\C_k|_\E)\times_{\R^s}(\C_l|_\mathcal{R})$.

 \begin{rk}
For $l=0$ we get extension of the distribution via the Cauchy characteristic, i.e. locally
$(M,\Delta)=(\bar M,\bar\Delta)\times(\R,\R)$.
 \end{rk}

This construction allows realizing all cases not prohibited by
Theorem \ref{Thm2} as distributions with infinite-dimensional symmetry algebras.

For rank 3 distributions all infinite primitive symmetry algebras
 % (non-primitive infinite cases include fiber products)
that occur are Lie transformations\footnote{In their natural representation
$k=0$ and $k=0,1$ in the next cases to be primitive,
but the algebra does not change with $k$ due to Lie-B\"acklund theorem.} of $J^k(\R,\R^2)$,
while for rank 4 there appear new real primitive algebras from the list of \cite{MT}:
the symmetries of the Cartan distribution of $J^k(\R,\R^3)$, $J^1(\R^2,\R)$ and $J^k(\CC,\CC)_\R$.

\subsection{Distributions of rank $n=3$}

Let us describe in more details the case of $\op{rank}(\Delta)=3$.
If the growth vector is $(3,1,..)$, then the bracket
$\La^2\g_{-1}\to\g_{-2}$ has a kernel $v\in\g_{-1}$, and this corresponds to a Cauchy characteristic
vector field, so that $(M,\Delta)\simeq(\bar M,\bar\Delta)\times(\R,\R)$.

Consider the growth $(3,2,..)$. Then $\Delta$ contains a rank 2 sub-dis\-trib\-ution $\Pi$, with
$\Pi_2\subset\Delta$. If $\Pi_2=\Delta$, then $\op{sym}(\Delta)=\op{sym}(\Pi)$. The growth vector of
$\Pi$ starts $(2,1,2,..)$ and so $\op{sym}(\Pi)$ is finite-dimensional.

Thus infinite algebras $\op{sym}(\Delta)$ correspond to integrable $\Pi$. In this case
the characteristic variety is real and $\op{Char}^\CC=\op{Char}$ equals the one
point set $\mathbb{P}(\Pi^\perp)$, where $\Pi^\perp$ is the annihilator of $\Pi\subset\g_{-1}$.

Moreover let $p\in\Pi^\perp\cap\g_{-1}^*\setminus0$. Then the kernel bundle over characteristic variety
$\mathcal{K}=\{q\in\g_{-1}:p\ot q\in\h_0\}$ is either 1- or 2-dimensional subspace of $\Pi$.

\smallskip

{\bf I:}  $\dim\mathcal{K}=1$. Here $\op{sym}(\Delta)$ is the symmetry of the flag
$(\mathcal{K},\Pi,\Delta)$. Moreover this triple extends (locally) to a complete flag
$\mathfrak{F}$ of subspaces $(\mathcal{K},\Pi,\Delta,[\mathcal{K},\Delta],\Delta_2,\dots)$
in $TM$, invariant under $\op{sym}(\Delta)$.

 \begin{rk}
For $n=2$ the infinite symmetry algebra also leaves invariant the complete flag
$\mathfrak{F}=(\langle\zeta\rangle,\Delta,\Delta_2,\dots)$, see Section \ref{SS81}.
 \end{rk}

Now $\mathcal{K}$ is the Cauchy characteristic space of the distribution $\Delta^\dag=[\mathcal{K},\Delta]$,
so we can pass to the (local) quotient $\bar M=M/\mathcal{K}$, $\bar\Delta=\Delta^\dag/\mathcal{K}$.
Any symmetry of $\Delta$ descends to a symmetry of $\bar\Delta$, and the latter contains a
preferred direction $\Pi/\mathcal{K}$. It is not necessarily the kernel direction $\bar{\mathcal{K}}$ of $\bar\Delta$, but it belongs to the corresponding 2-distribution $\bar\Pi$.

The passage $\Delta\mapsto\bar\Delta$ is de-prolongation of rank 3 distributions in the same sense as
in Section \ref{SS81}, namely this operation is inverse to the prolongation defined as follows
(this applies only for $\op{Char}(\Delta)\ne\emptyset$, so that
$\Pi\subset\Delta$ is integrable etc).

\underline{Prolongation I}$_a$. This works in the case $\bar\Delta$ has a preferred direction
$\ell=\langle Y\rangle$ different from $\bar{\mathcal{K}}$
[if not, then the Lie algebra $\op{sym}(\bar\Delta)$ shrinks to the stabilizer].
Namely $M$ is the space of all 2-planes $P\subset\bar\Delta$ through $\ell$ transversal to $\bar{\mathcal{K}}$
(this is the affine version; compact version is without transversality).
If $\pi:M\to\bar M$ is the projection, then $\Delta=\pi_*^{-1}(P)$.

Locally if $\bar\Delta=\langle X,Y,Z\rangle$ with $X\in\bar{\mathcal{K}}$ and $Y\in\ell\subset\bar\Pi$,
then $M=\bar M\times\R(t)$, $\Delta=\langle Y,\p_t,Z+tX\rangle$.
In these notations $\Pi=\langle Y,\p_t\rangle$, $\p_t\in\mathcal{K}$.
This is inverse to the above de-prolongation when $\dim\mathcal{K}=1$.

\underline{Prolongation I}$_b$. $\bar{\mathcal{K}}$ is clearly a preferred direction,
so we can prolong along it: $M$ is the space of all 2-planes $P\subset\bar\Delta$ through
$\bar{\mathcal{K}}$ different from $\bar{\Pi}$
(the affine version; compact version - no conditions).

Again $\Delta=\pi_*^{-1}(P)$ and locally $\Delta=\langle X,\p_t,Z+tY\rangle$. Here
$\Pi=\langle X,\p_t\rangle=\mathcal{K}$, so this is not inverse to the above de-prolongation.

\smallskip

{\bf II:}  $\dim\mathcal{K}=2$. Here $\Pi=\mathcal{K}$ is the space of Cauchy characteristics
for rank 5 distribution $\Delta_2$, and we can reduce $(M,\Delta)$ to a rank 3 distribution
$\bar\Delta=\Delta_2/\Pi$ on a lower-dimensional manifold $\bar M=M/\Pi$.
This is de-prolongation in the same sense as in Section \ref{SS81}, while the prolongation
($\Pi$ integrable) mimics the Spencer prolongation on $J^k(\R,\R^2)$.

\underline{Prolongation II}. Here $M$ is the space of lines $L\subset\Delta$ transversal to
$\Pi$ (the affine version; compact version - no restrictions), $\pi:M\to\bar M$ the natural projection,
then $\Delta=\pi_*^{-1}(L)$.

Locally if $\bar\Delta=\langle X,Y,Z\rangle$ with $\Pi=\langle X,Y\rangle$, then $M=\bar M\times\R^2(u,v)$
and $\Delta=\langle \p_u,\p_v,Z+uX+vY\rangle$. Clearly $\mathcal{K}=\langle\p_u,\p_v\rangle$.

\smallskip

Notice that de-prolongations in both cases induce the injective map
$\pi_*:\op{sym}(\Delta)\hookrightarrow\op{sym}(\bar\Delta)$
[prolongations I$_a$ can shrink the symmetry algebra, though in a controllable way].

Thus de-prolongations can be applied until the distribution gets a Cauchy characteristic
or its length (degree of non-holonomy) becomes $\kappa=1$.
Thus we obtain the following normal form [for $n=3$ infinitely-symmetric models can have
functional moduli, but these are gone for most symmetric models in the sense of functional dimension].

 \begin{theorem}
A regular germ of a rank 3 non-holonomic distribution $\Delta$ has infinite-dimensional symmetry
algebra if and only if $\Delta$ is obtained from the product $\Pi\times\R$, where
$\Pi$ is some rank 2 distribution, by the operations of prolongations (of type I or II).
 \end{theorem}

Here we include the case when $\Pi$ is the trivial distribution $(\R^2,\R^2)$, in which case prolongation II
of $\Pi\times\R$ yields $J^k(\R,\R^2)$.

\smallskip

 {\bf Example I.} Consider the Monge equation $\E_{m,n}=\{y_m=z_n^2\}$ as a submanifold
in the mixed jets (subscripts count the derivatives by $x$)
 $$
J^{m,n}(\R,\R^2)\simeq\R^{m+n+3}(x,y,y_1,\dots,y_m,z,z_1,\dots,z_n)
 $$
and let $M=\E_{m,n}\times_\R J^l(\R,\R)\subset J^{m,n,l}(\R,\R^3)$
with the natural rank 3 distribution
$\Delta=\mathcal{C}_{m,n}\times_\R\mathcal{C}_l=\langle\D_x,\p_{z_n},\p_{w_l}\rangle$,
 $$
\D_x=\p_x+y_1\p_y+{}\dots{}+z_n^2\p_{y_{m-1}}+z_1\p_z+{}\dots{}
+z_n\p_{z_{n-1}}+w_1\p_w+{}\dots{}+w_n\p_{w_{l-1}}.
 $$
De-prolongation I results in
$\bar M=\R^{m+n+l+2}(x,\{y_i\}_0^{m-1},\{z_i\}_0^n,\{w_i\}_0^{l-1})$
with $\bar\Delta=\langle\D_x,\p_{z_n},\p_{w_{l-1}}\rangle$.
Further de-prolongations I bring the distribution to
$\hat\Delta=\langle\D_x,\p_{z_n},\p_w\rangle=\mathcal{C}_{m,n}\times\R$ on
$\E_{m,n}\times\R(w)$.

The general symmetry of this distribution is $X_0=Y_\E+f\cdot\p_w$, where $Y_\E$ is
the general symmetry of $\E_{m,n}$, $m\le n$ (it depends on $2n+5$ parameters for $m=1$, $n>2$
or $2n+4$ parameters for $m>1$, see \cite{AK}) and $f\in C^\infty(\E_{m,n}\times\R)$.

In order for this symmetry to allow prolongation I$_a$, it shall preserve the vertical
line $\mathfrak{v}=\langle\p_{z_n}\rangle$. This is so for the first component $Y_\E$
\cite{AK}, and it holds for $f\cdot\p_w$ iff $f_{z_n}=0$. The prolongation is
$X_1=Y_\E+f\cdot\p_w+\D_x(f)\cdot\p_{w_1}$. Now this preserves $\mathfrak{v}$ iff
$f_{z_{n-1}}=0$, $f_{y_{m-1}}=0$, and the prolongation is
$X_2=Y_\E+f\cdot\p_w+\D_x(f)\cdot\p_{w_1}+\D_x^2(f)\cdot\p_{w_2}$.

Continuing in this way we obtain that the symmetries of $(M,\Delta)$ are
$X_l=Y_\E+\sum_{k\le l}\D^k_x(f)\cdot\p_{w_k}$ and $f=f(x,\{y_i\}_0^{m-l},\{z_i\}_0^{n-l},w)$
is any smooth function of $2+\max\{0,m-l+1\}+\max\{0,n-l+1\}$ arguments.

\smallskip

 {\bf Example II.} For the mixed jets $J^{m,n}(\R,\R^2)=J^m(\R,\R)\times_\R J^n(\R,\R)$ with
canonical Cartan distribution $\mathcal{C}_{m,n}$ of rank 3, de-prolongation II reduces it to
$J^{0,n-m}(\R,\R^2)=\R\times J^{n-m}(\R,\R)$, $m\le n$. The infinite algebra of symmetries
for this model is described in \cite{AK}.

\smallskip

Some other results on distributions with infinitely many symmetries
(on graded nilpotent Lie groups) can be found in \cite{DR,OR}.

%%%%%%%%%%%%%%%
\appendix
\section{Estimate of the solution space size}\label{A1}

In the theory of formal integrability \cite{Go,Sp} a PDE (system) is
considered as a submanifold $\E_k$ in jets $J^k\pi$ of a bundle $\pi:E\to M$
(if different orders are involved, $\E$ is a collection of submanifolds \cite{KL}),
and its prolongations defined geometrically are submanifolds, possibly with singularities,
in higher jets $\E_i\subset J^i\pi$.

The (tangent to the) fiber of the projection $\pi_{i,i-1}:\E_i\to\E_{i-1}$
is the symbol space $g_i\subset S^iT^*\ot N$. Then the following estimate on the
size of the solution space holds formally (i.e. the solutions are viewed as formal
series at $x$ and the r.h.s. is evaluated at the same point):
 \begin{equation}\label{est-g}
\dim\op{Sol}(\E)\le\sum_0^\infty\dim g_i.
 \end{equation}
Of course, this inequality is meaningful only in the case the r.h.s. is finite
(such $\E$ are said to have finite type at $x$).

If the equation $\E$ is formally integrable and regular ($\E_i$ are smooth and
dimensions of $g_i$ are constant), the same estimate holds locally and globally
(but in the latter case the l.h.s. can shrink by other reasons).

But often for a regular but formally integrable overdetermined PDE $\E$
the result of prolongation-projection $\bar\E$ (which is formally integrable
in regular points by Cartan-Kuranishi theorem \cite{Ku}) is not everywhere
regular (e.g. for Lie equations). We would like to prove inequality
(\ref{est-g}) in the presence of singularities under some mild assumptions.

To be more precise, assume that the original equation $\E_k$ is regular
(we restrict to pure order systems, but the result can be generalized),
consider the nested sequence $\E_i^l=\pi_{l,i}(\E_k^{(l-k)})$ and let $\bar\E_i=\cap_l\E_i^l$.
Then the co-filtration $\bar\E_i$ represents a formally integrable equation
away from singularities.

We restrict to linear equations $\E$, as our main application here -- the Lie equation
for symmetries $\frak{Lie}(K)$ -- is linear.

Our main result states that if in finite number of steps of prolongation-projection
we get finite type at all points (including singularities),
then we conclude (\ref{est-g}) even before arriving to formally integrable equation.

Denote by $\bar g_i$ the fiber of the projection $\pi_{i,i-1}:\bar\E_i\to\bar\E_{i-1}$
(can be empty). For linear equations the symbol $\bar g_i$ depends on $x\in M$.

 \begin{theorem}\label{T-app}
Let $\E_k\subset J^k(\pi)$ be a linear equation in a bundle over a connected manifold $M$.
The space $\op{Sol}(\E)$ of global solutions over $M$ is a linear space, and its
dimension is well-defined.
Suppose that there exists a number $l$ such that the projection
$\pi_{i,i-1}:\E_i\to\E_{i-1}$ is injective for $i=l$ (and hence for $i\ge l$), where
$\E_i=\E_k^{(i-k)}$. Then
 \begin{equation}\label{est-g2}
\dim\op{Sol}(\E)\le\inf_M\sum_0^\infty\dim\bar g_i(x).
 \end{equation}
 \end{theorem}

 \begin{proof}
The result follows from the following claim (with the imposed assumptions):
for every point $x\in M$ the map $j_x^\infty:\op{Sol}(\E)\to J_x^\infty\pi$ is
injective. Indeed, in this case $\op{Sol}(\E)$ is a subset of $\bar\E_x$, which has dimension
bounded by $\sum_0^\infty\dim\bar g_i(x)$. What the claim states is that there
exists no solutions flat at the given point.

To prove the claim consider a relatively compact connected neighborhood $U''\subset\E_{l+1}$
and let $U'=\pi_{l+1,l}(U'')\subset\E_l$, $U=\pi_{l,l-1}(U')\subset\E_{l-1}$.
For any $x_l\in U'$ and $x_{l-1}=\pi_{l,l-1}(x_l)\in U$ we have the horizontal $n$-plane
$\Pi(x_{l-1})=d\pi_{l,l-1}(\C(x_l))\subset\C(x_{l-1})\subset T_{x_{l-1}}J^{l-1}\pi$,
where $\C(\cdot)$ is the canonical Cartan distribution in jets (\cite{KLV}).

If $U$ is sufficiently small, then $\pi_{l,l-1}:U'\to U$ is a homeomorphism. Indeed, it
is bijective (the equation is linear, the fiber $g_l$ of $\pi_{l,l-1}$ is zero) and
is a continuous map of a compact into a Hausdorff space. Moreover
the slope of $\Pi(x_{l-1})$ (in canonical coordinates) is determined by $x_l$.
Variation of slopes is governed by $x_{l+1}=\pi_{l+1,l}^{-1}(x_l)\cap U''$,
and $\pi_{l+1,l}:U''\to U'$ is a homeomorphism as well.
Consequently this distribution of planes has Lipschitz dependence on
the point $x_{l-1}\in\E_{l-1}$ (and is $C^1$ in regular points).

Suppose there are two local solutions of $\E$ with infinite tangency at $x\in M$.
Let $x_{l-1}$ be their $(l-1)$-jet at $x$, and denote by $\Upsilon,\Upsilon'\subset\E_{l-1}$
the corresponding integral manifolds (jet-lifts) through $x_{l-1}$.

Consider a curve $\sigma\subset M$ through $x$ and denote $M_\z=\pi_{l-1}^{-1}(\z)\cap U$.
Then $\gamma_\z=\Upsilon\cap M_\z$ and $\gamma'_\z=\Upsilon'\cap M_\z$ are integral curves of
the Lipschitz line distribution $l_\z=\Pi\cap TM_\z$.

As in the usual uniqueness theorem for ODEs this implies our claim. Actually
if $\gamma_\z(t),\gamma'_\z(t)$ are parametrized integral curves, $t\in I$ and
$\t=|I|$ is the length of the interval, then for $C^0$-max norms
(in any coordinates) we have from the integral version of the ODE $\dot\gamma=l_\z\circ\gamma$:
 $$
\|\gamma_\z-\gamma_\z'\|\le c\t\,\|\gamma_\z-\gamma_\z'\|
 $$
with some constant $c$ depending only on $\Pi$ and $\z$. Shrinking $I$ if necessary we obtain $c\t<1$
and so locally $\gamma_\z=\gamma_\z'$. As $\z$ is arbitrary, we get locally $\Upsilon=\Upsilon'$.
As $U$ is connected, this implies coincidence over the whole neighborhood $\pi_l(U)$
and further over the whole domain of the solutions (ultimately $M$).
 \end{proof}

A version of this theorem holds for non-linear $\E$ (then $\dim$ in the l.h.s. means the
Hausdorff dimension, the infimum in the r.h.s. is taken over the equation $\bar\E$ and
the symbols $\bar g_i$ are the linear envelopes of the $\pi_{i,i-1}$-vertical tangent cone
to $\bar\E_i$ at the singular point $x_i$), but more care shall be taken about regularity of the original equation $\E$.

Let us notice that the statement obviously holds true in the case, when the singularities in the
differential closure $\bar\E$ of $\E$ have $\op{codim}\ge2$.

\smallskip

In the opposite case we have the following counter-example (thanks to Valentin Lychagin).

{\bf Example.} Consider the ODE $\E:\ x^3\,y'(x)=y$. Then $g_0=\R$ and $g_i=0$ for $i>0$ and
all points $x\ne0$. Near all such points the equation is regular and the solution space
is one-dimensional. However the space of global solutions $\op{Sol}(\E)$ is two-dimensional
with the basis $\chi(x)e^{-1/x^2},\chi(-x)e^{-1/x^2}$,
where $\chi(x)$ is the Heaviside function.

Notice however that $\E$ is neither regular, nor of finite type at 0, and the jet-lift
does not resolve this singularity. In addition, the codimension of the singularity is 1.

 \begin{rk}
Thus a generalization of (\ref{est-g2}) could fail if one wants to disregard the singularities
and change the r.h.s. to $\op{ess.sup}_M\sum_0^\infty\dim\bar g_i(x)$.
 \end{rk}

A modification of the above example yields an equation $\E$ having finite type at all
regular points but with infinite-dimensional $\op{Sol}(\E)$.

%%%%%%%%%%%%%%%%%%%%%%%%%%%%%%%%%%%%%%%%%%%%%%%%%%%%%%%%%%%%%%%%%%%%%%%%%%%%

\end{document}